\documentclass[10pt]{article}
\usepackage{theorem,float}
\usepackage{amsmath,amssymb}
\usepackage{mathrsfs}
\usepackage[pdftex]{graphicx}
\restylefloat{figure}
\newcommand{\Proof}{{\bf Proof:\quad}}
\newcommand{\EndProof}{\hspace*{\fill} \rule{6pt}{6pt}\\[1\baselineskip]}
\newtheorem{theorem}{\bf Theorem}

\newtheorem{lemma}[theorem]{\bf Lemma}
\newtheorem{remark}[theorem]{Remark}
\newtheorem{corollary}[theorem]{Corollary}
\newcommand{\Lemma}{\begin{lemma}\quad}
\newcommand{\EndDefinition}{\end{defi}}
\newcommand{\EndLemma}{\end{lemma}}

\newcommand{\BR}{\mathbb{R}}

\title{Approximation Properties of 
a Gradient Recovery Operator Using a Biorthogonal System}
\author{Bishnu P. Lamichhane\footnotemark[1]\;  and Adam McNeilly\thanks{School of Mathematical and Physical Sciences, University of Newcastle, 
Callaghan, NSW 2308,  {\tt bishnu.lamichhane@newcastle.edu.au, adam.mcneilly@uon.edu.au
}}}
\begin{document}
\maketitle
\begin{abstract}

A gradient recovery operator based on projecting 
the discrete gradient onto the standard finite element space is considered.
We use an oblique projection, where the test and 
trial spaces are different, and the bases of these two spaces form
a biorthogonal system.  Biorthogonality allows efficient computation of 
the recovery operator.
We analyse the 
approximation properties of 
the gradient recovery operator. 
\end{abstract}

\section{Introduction}
The gradient reconstruction is a popular technique to develop a 
reliable a posteriori error estimators for approximating the solution 
of partial differential equations using adaptive finite element methods 
\cite{ZZ92a,ZZ92b,AO00,XZ04,Lon06,CW10}.
Recently we have presented a gradient reconstruction operator based on an oblique projection 
\cite{Lam10a}. The oblique projection operator is constructed by 
using a biorthogonal system. In fact, for the linear finite element in simplicial meshes, 
this approach reproduces the so-called gradient reconstruction scheme by 
the weighted averaging \cite{Goo94,XZ04,CW10}. We proved that 
the approximation property of the recovered gradient for any finite element space 
is similar to the one obtained by using the orthogonal projection with respect 
to $L^2$-norm \cite{XZ04,Lon06}. 
 In this article, we aim at analysing the approximation property 
of the recovered gradient in one dimension using the oblique projection.
This construction is quite useful in extending the 
weighted average gradient recovery of linear finite elements \cite{XZ04}
to quadrilaterals and hexahedras. 

Let $\Omega =(\alpha,\beta)$ with $\alpha,\beta \in \BR$ and $\alpha<\beta$. Let
$\Delta=\left\{ \alpha=x_{0}<x_{1}<...<x_{n}=\beta\right\} $ be a partition 
of the interval $\Omega$. We define the interior of the grid, denoted
$int\left(\Delta\right)$, as
\[
int\left(\Delta\right)=\left\{ x_{i}\in\Delta\,:\,1\leq i\leq n-1\right\} .
\]
We also define the set of intervals in the partition $\Delta$ as 
$ \left\{ I_{i}\right\} _{i=0}^{n-1}$, where
$I_{i}=[x_{i},\, x_{i+1})$. 
Two sets $A_n$ and $B_n$ of indices are also defined as 
$A_{n}=\left\{ i\in\mathbb{N}\,:\,0\leq i\leq n\right\} $ and
$B_{n}=A_{n}\setminus\left\{ 0,\, n\right\} $, respectively. 
A piecewise linear interpolant of a continuous function \emph{$u$
} is written as $I_{h}u\in V_{h}$ with
\[
I_{h}u(x)=\sum_{i=0}^{n}u(x_{i})\phi_{i}(x),
\]
where $ \phi_i$ is the standard hat function associated with the 
point $x_i$, $0 \leq i \leq n$. We define a discrete space, 
\[ V_{h}=\mathrm{span}\left\{
  \phi_{0},\,...\,,\phi_{n}\right\} \subset H^{1}\left(\Omega\right).\]
The linear interpolant of $u \in H^1(\Omega)$ is the 
continuous function defined by 
$I_h u =\sum_{i=0}^{n}u(x_{i})\phi_{i}$.
However, if we compute the derivative of this interpolant 
$I_hu$, the resulting function will not be continuous. 
To make the derivative continuous 
we  project the 
derivative of the interpolant, $\frac{\partial u_{h}}{\partial x}=\sum_{i=0}^{n}u(x_{i})\,\frac{\partial\phi_{i}}{\partial x}$,
onto the discrete space $V_{h}$.   
There are two different types of projection. One is an orthogonal projection and the 
other is an oblique projection. 
The orthogonal projection operator, $P_{h}$, that projects $\frac{\partial u_{h}}{\partial x}$
onto $V_{h}$ is to find a $g_{h}=P_{h}\frac{\partial u_{h}}{\partial x}\in V_{h}$
that satisfies:
\begin{equation}
\int_{\Omega}g_{h}\phi_{j}\, dx=\int_{\Omega}\frac{\partial u_{h}}{\partial x}\phi_{j}\, dx.\label{eq: 2.1}
\end{equation}

Since $g_{h}\in V_{h}$, we can represent it as an $\left(n+1\right)$-dimensional
vector:

\[
\vec{g}=\left(\begin{array}{c}
g_{0}\\
\vdots\\
g_{n}
\end{array}\right)\,\mathrm{with}\; g_{h}=\sum_{i=0}^{n}g_{i}\phi_{i}.
\]

Now the requirement given in equation (\ref{eq: 2.1}) is equivalent
to a linear system: $M\vec{g}=\vec{f}$, where $M$ is a mass matrix,
and

\[
\vec{f}=\left(\begin{array}{c}
f_{0}\\
\vdots\\
f_{n}
\end{array}\right)\,\mathrm{with}\: f_{j}=\int_{\Omega}\frac{\partial u_{h}}{\partial x}\phi_{j}\, dx.
\]
 Here the mass matrix $M$ is tridiagonal. We can reduce computation time
greatly if we have a diagonal mass matrix. This can be done if we
use  a suitable oblique
projection instead of an orthogonal projection. We consider the projection
\[ Q_{h}:L^2(\Omega)\rightarrow V_h,\] which is defined as the problem of finding 
  $g_{h}=Q_{h}\frac{du}{dx}\in V_{h}$ 
such that
\[
\int_{\Omega}g_{h}\lambda_{h}\, dx=\int_{\Omega}\frac{du_{h}}{dx}\lambda_{h}\, dx,\qquad\lambda_{h}\in M_{h},
\]
where $M_{h}$ is another piecewise polynomial space, not orthogonal
to $V_{h}$, with $\dim\left(M_{h}\right)=\dim\left(V_{h}\right)$,
see \cite{Woh01}. In fact, the projection operator $Q_h$ is well-defined due to 
the following stability condition.
There is a constant $\beta>0$ independent of the mesh-size $h$  such that \cite{Lam10a,Lam11b}
\begin{eqnarray*}
\|\phi_h\|_{L^2(\Omega)} \leq \beta \sup_{\mu_h \in M_h \backslash\{0\}} 
\frac{\int_{\Omega} \mu_h\phi_h\,dx} {\|\mu_h\|_{L^2(\Omega)}},
\quad \phi_h \in V_h.
\end{eqnarray*}
In order to achieve that the mass matrix $M$ is diagonal 
we need to define a new set of basis functions
for $M_{h}$, $\left\{ \lambda_{0},\,...\,,\lambda_{n}\right\} $,
that are biorthogonal to the standard hat basis function we used previously.
This biorthogonality relation is defined as:
\begin{equation}
\int_{\Omega}\lambda_{i}\phi_{j}dx=c_{j}\delta_{ij},\, c_{j}\not=0,\,1\leq i,\, j\leq n,\label{eq: 2.2}
\end{equation}
where $\delta_{ij}$ is the Kronecker delta function:
\[
\delta_{ij}=\begin{cases}
1, & \mathrm{if}\, i=j\\
0, & \mathrm{otherwise}
\end{cases},
\]
and $c_{j}$ is a positive scaling factor. The basis functions for $M_h$ are simply given by 
\begin{align*}
\lambda_{0}(x) & =\begin{cases}
\frac{2(x-x_{1})+(x-x_{0})}{x_{0}-x_{1}}, & x_{0}\leq x\leq x_{1}\\
0\,, & \mathrm{otherwise}
\end{cases}\\
\lambda_{n}(x) & =\begin{cases}
\frac{2(x-x_{n-1})+(x-x_{n})}{x_{n}-x_{n-1}}, & x_{n-1}\leq x\leq x_{n}\\
0\,, & \mathrm{otherwise}
\end{cases},
\end{align*}
and for $1\leq i \leq n-1$
\begin{align*}
\lambda_{i}(x) & =\begin{cases}
\frac{2(x-x_{i-1})+(x-x_{i})}{x_{i}-x_{i-1}}, & x_{i-1}\leq x\leq x_{i}\\
\frac{2(x-x_{i+1})+(x-x_{i})}{x_{i}-x_{i+1}}, & x_{i}\leq x\leq x_{i+1}\\
0\,, & \mathrm{otherwise}
\end{cases}.
\end{align*}

\begin{figure}[H]
\begin{centering}
\includegraphics[scale=0.3]{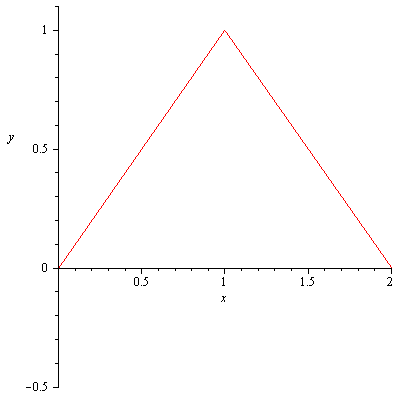}\includegraphics[scale=0.3]{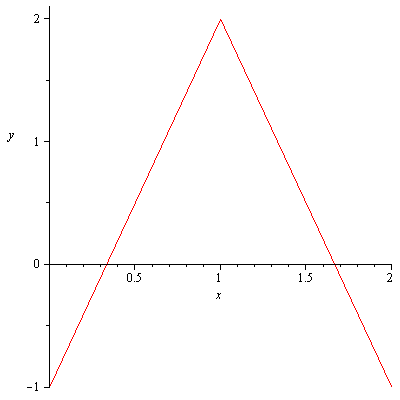}
\par\end{centering}

\caption{The hat basis function (left) and biorthogonal basis function
  (right)  with stepsize $h=1$.}
\end{figure}

By using an oblique projection $Q_{h}$ the mass matrix will be diagonal.
We let the diagonal mass matrix be D, so that our system is $D\vec{g}=\vec{f}$.
The values $g_{i}$ are our estimates of the gradient of $u$ at the
point $x_{i}$. So, we estimate the gradient by finding $\vec{g}=D^{-1}\vec{f}$,
where\[
g_{i}=\frac{\int_{x_{0}}^{x_{n}}\frac{du_{h}}{dx}\lambda_{i}dx}{\int_{x_{0}}^{x_{n}}\phi_{i}\lambda_{i}dx}.
\]
We want to calculate the error in this approximation, and find out
when $g_{i}$ approximates $u'(x_{i})$ exactly for each $x_{i}\in\Delta$.
As in \cite{Zha96,Zha00} we want to see if
 $g_{i}$ approximates $u'(x_{i})$ exactly when $u$ is a
quadratic polynomial.

\section{Superconvergence}
\begin{theorem}\label{th0}
Let $u\in C^0(\Omega)$. Then we have 
\[ g_i = \frac{u(x_{i+1}) - u(x_{i-1})}{x_{i+1}-x_{i-1}}, \quad i \in B_n, \]
and 
\[ g_0 =\frac{u(x_{1}) - u(x_0)}{x_{1}-x_{0}}, 
\quad 
g_n = \frac{u(x_{n}) - u(x_{n-1})}{x_{n}-x_{n-1}}.\]
\end{theorem}
\Proof
We note that 
\[
\lambda_{i}(x)=\begin{cases}
\frac{2(x-x_{i-1})+(x-x_{i})}{x_{i}-x_{i-1}}, & x_{i-1}\leq x\leq x_{i}\\
\frac{2(x-x_{i+1})+(x-x_{i})}{x_{i}-x_{i+1}} & x_{i}\leq x\leq x_{i+1}\\
0\,, & \mathrm{otherwise}
\end{cases},\:\mathrm{for\, all}\: i\in B_{n}.
\]
Now, we calculate $g_{i}$ for $i\in B_{n}$: 
\[
g_{i}=\frac{\int_{x_{0}}^{x_{n}}\frac{du_{h}}{dx}\lambda_{i}\, dx}{\int_{x_{0}}^{x_{n}}\phi_{i}\lambda_{i}\, dx},
\]
where 
\[
\begin{aligned} \int_{x_{0}}^{x_{n}}\phi_{i}\lambda_{i}dx & =\int_{x_{i-1}}^{x_{i}}\left(\frac{x-x_{i-1}}{x_{i}-x_{i-1}}\right)\left(\frac{2\left(x-x_{i-1}\right)+\left(x-x_{i}\right)}{x_{i}-x_{i-1}}\right)dx\\
 & +\int_{x_{i}}^{x_{i+1}}\left(\frac{x-x_{i+1}}{x_{i}-x_{i+1}}\right)\left(\frac{2\left(x-x_{i+1}\right)+\left(x-x_{i}\right)}{x_{i}-x_{i+1}}\right)dx\\
 & =-\frac{1}{2}\left(x_{i-1}-x_{i+1}\right),
\end{aligned}
\]
and
\[
\begin{aligned} \qquad\int_{x_{0}}^{x_{n}}\frac{du_{h}}{dx}\lambda_{i}\, dx & =\sum_{j=0}^{n}u(x_{j})\int_{x_{0}}^{x_{n}}\frac{d\phi_{j}}{dx}\lambda_{i}\, dx\\
 & =u(x_{i-1})\int_{x_{i-1}}^{x_{i}}\frac{d\phi_{i-1}}{dx}\lambda_{i}\, dx+u(x_{i})\left(\int_{x_{i-1}}^{x_{i}}\frac{d\phi_{i}}{dx}\lambda_{i}\, dx+\int_{x_{i}}^{x_{i+1}}\frac{d\phi_{i}}{dx}\lambda_{i}\, dx\right)\\
 & +u(x_{i+1})\int_{x_{i}}^{x_{i+1}}\frac{d\phi_{i+1}}{dx}\lambda_{i}\, dx\\
 & \mathrm{(since}\:\phi_{j}\:\mathrm{and}\:\lambda_{i}\:\mathrm{overlap\, only\, when}\, j\in\left\{ i-1,\, i,\, i+1\right\} \mathrm{)}\\
\mathrm{Therefore,\quad}
 g_i &= \frac{u(x_{i+1}) - u(x_{i-1})}{x_{i+1}-x_{i-1}}.
\end{aligned}
\]Now we look at the end-points.  We note that 
\[
g_{0}=\frac{\int_{x_{0}}^{x_{n}}\frac{du_{h}}{dx}\lambda_{0}\, dx}{\int_{x_{0}}^{x_{n}}\phi_{0}\lambda_{0}\, dx},\quad 
\text{and}\quad 
g_{n}=\frac{\int_{x_{0}}^{x_{n}}\frac{du_{h}}{dx}\lambda_{n}\, dx}{\int_{x_{0}}^{x_{n}}\phi_{n}\lambda_{n}\, dx}.
\]
Computing as before we get 
\[ g_0 =\frac{u(x_{1}) - u(x_0)}{x_{1}-x_{0}}, 
\quad 
g_n = \frac{u(x_{n}) - u(x_{n-1})}{x_{n}-x_{n-1}}.\]
\EndProof
We have the following super-convergence in $L^2$-norm. This is 
proved as in \cite{LZ99,XZ04}.
\begin{theorem}
 Let $h_i = x_{i+1} - x_{i}$ for $0 \leq i \leq n-1$, $h = \max_{0 \leq i \leq n-1}h_i$, 
 and $\Omega_0 = (x_1,x_{n-1})$. 
 If the point distribution satisfies $|h_{i+1} - h_i| =O(h^2)$ for $0 \leq i \leq n-1$. 
Then we have the estimate  
 \[ \left\|\frac{du}{dx} - Q_h\frac{du_h}{dx}\right\|_{L^2(\Omega_0)} \leq 
 h^2 \|u\|_{W^{3,\infty}(\Omega_0)},\quad 
 u \in W^{3\infty}(\Omega).\] 
\end{theorem}
For the tensor product meshes in two or three dimensions satisfying the 
above mesh condition this theorem has an easy extension. 

\subsection{Application to  quadratic functions}
\begin{corollary}\label{th1}
Let $u\in P_{2}\left(\mathbb{R}\right)$. Then $g_{i}$ reproduces
$u'\left(\tilde{x}_{i}\right)$ exactly for all $x_{i}\in\Delta$,
where:
\begin{eqnarray*}
\tilde{x}_{i}=\begin{cases}
\frac{x_{0}+x_{1}}{2}, & i=0\\
\frac{x_{i-1}+x_{i+1}}{2}, & 1\leq i\leq n-1\\
\frac{x_{n-1}+x_{n}}{2}, & i=n
\end{cases},\quad\text{and}\quad 
g_{i}=\frac{\int_{x_{0}}^{x_{n}}\frac{du_{h}}{dx}\,\lambda_{i}\, dx}{\int_{x_{0}}^{x_{n}}\phi_{i}\,\lambda_{i}\, dx}.
\end{eqnarray*}
\end{corollary}

\begin{figure}[H]
\begin{centering}
\includegraphics[scale=0.8]{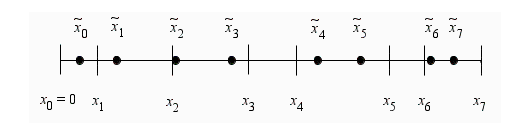}
\end{centering}
\caption{A non-uniform grid with 8 nodes (vertical lines). The points $\tilde{x}_{i}$ (dots)
are also shown.}\end{figure}
\Proof
We use the result of the previous theorem to get 
\[
\begin{aligned} 
g_{i} & =\frac{\left(ax_{i-1}^{2}+bx_{i-1}-ax_{i+1}^{2}-bx_{i+1}\right)}{\left(x_{i-1}-x_{i+1}\right)}\\
 & =\frac{a\left(x_{i-1}^{2}-x_{i+1}^{2}\right)+b\left(x_{i-1}-x_{i+1}\right)}{x_{i-1}-x_{i+1}}\\
 & =\frac{a\left(x_{i-1}-x_{i+1}\right)\left(x_{i-1}+x_{i+1}\right)+b\left(x_{i-1}-x_{i+1}\right)}{x_{i-1}-x_{i+1}}\\
 & =a\left(x_{i-1}+x_{i+1}\right)+b.
\end{aligned}
\]
On the other hand,  \[
\begin{aligned}u'\left(\tilde{x}_{i}\right) & =2a\left(\frac{x_{i-1}+x_{i+1}}{2}\right)+b\\
 & =a\left(x_{i-1}+x_{i+1}\right)+b\\
 & =g_{i}.
\end{aligned}
\]
So, $g_{i}$ reproduces $u'\left(\tilde{x}_{i}\right)$ exactly for
$i\in B_{n}$. Now for $ i =0$ and $i=n$, we have 
\[
\begin{aligned}
g_{0} & =\frac{\frac{1}{2}\left(ax_{1}^{2}+bx_{1}-ax_{0}^{2}-bx_{0}\right)}{\frac{1}{2}\left(x_{1}-x_{0}\right)}=\frac{a\left(x_{1}^{2}-x_{0}^{2}\right)+b\left(x_{1}-x_{0}\right)}{x_{1}-x_{0}} \\ & =\frac{a\left(x_{1}-x_{0}\right)\left(x_{1}+x_{0}\right)+b\left(x_{1}-x_{0}\right)}{x_{1}-x_{0}} =a\left(x_{1}+x_{0}\right)+b,\\
 \quad\text{and}\quad 
g_{n} & =\frac{\frac{1}{2}\left(ax_{n-1}^{2}+bx_{n-1}-ax_{n}^{2}-bx_{n}\right)}{\frac{1}{2}\left(x_{n-1}-x_{n}\right)}
  =\frac{a\left(x_{n-1}^{2}-x_{n}^{2}\right)+b\left(x_{n-1}-x_{n}\right)}{x_{n-1}-x_{n}}\\
 & =\frac{a\left(x_{n-1}-x_{n}\right)\left(x_{n-1}+x_{n}\right)+b\left(x_{n-1}-x_{n}\right)}{x_{n-1}-x_{n}}=a\left(x_{n-1}+x_{n}\right)+b.
\end{aligned}
\]
Since 
\[u'\left(\tilde{x}_{0}\right)  =a\left(x_{0}+x_{1}\right)+b,\quad\text{and}\quad 
u'\left(\tilde{x}_{n}\right)  =a\left(x_{n-1}+x_{n}\right)+b,
\]
we have 
$g_0$ and $g_n$ reproduce $u'\left(\tilde{x}_{0}\right)$  and 
$u'\left(\tilde{x}_{n}\right)$, respectively, exactly. 
\EndProof

\begin{remark}[Uniform Grid]
Let $\Delta=\left\{ \alpha=x_{0}<x_{1}<...<x_{n}=\beta\right\} $ be a uniform 
grid on the interval $\Omega$ so that 
$x_{i}-x_{i-1}=h,\;\forall\, i\in A_{n}\setminus\left\{ 0\right\} $,
where $h$ is some constant, called the \emph{stepsize}.
We note that if our grid is uniform, 
then $x_{i}=\tilde{x}_{i}\;\forall x_{i}\in int\left(\Delta\right)$.
So, our gradient recovery operator will reproduce the exact gradient
of any quadratic function on the interior of a uniform grid. We cannot
recover the gradients at the endpoints exactly, however, since $x_{0}\not=\tilde{x}_{0}$
and $x_{n}\not=\tilde{x}_{n}$. 
\end{remark}

\begin{corollary}
Let $u\in P_{2}(\mathbb{R})$ with $u(x)=ax^{2}+bx+c$, and let the
grid be uniform with stepsize $h$ and $x_{0}=0$. Then $\left|g_{i}-u'(x_{i})\right|=ah$
for $i=0,\, n$ (i.e. for the endpoints of the grid).
\end{corollary}
\Proof
We will start with the case where $i=0$ (i.e. the left endpoint).
We know from Theorem  \ref{th1}  that $g_{0}=a\left(x_{1}-x_{0}\right)+b$.
Since our grid is uniform with stepsize $h$, this simplifies to $g_{0}=ah+b$. 
$u'\left(x_{0}\right)=2ax_{0}+b=b$, since $x_{0}=0$.
Therefore
\begin{alignat*}{1}
\left|g_{0}-u'\left(x_{0}\right)\right| & =\left|\left(ah+b\right)-b\right|\\
 & =ah.
\end{alignat*}
The case for $i=n$ (i.e. the right endpoint) is proven similarly.
\EndProof

For a non-uniform grid, we cannot simplify our approximations using
the stepsize $h$, since the spacing between each adjacent node is
not always equal. We did not make any assumption about the uniformity of the grid 
in Theorem \ref{th1}. Thus  $|g_i -u'(x_i)|$, $ i \in A_n$, is not zero for a non-uniform grid 
This is estimated in the following corollary.

\begin{corollary}
Let $u\left(x\right)=ax^{2}+bx+c$.
Then,  \begin{alignat*}{1}
\left|g_{i}-u'\left(x_{i}\right)\right| & =a\left(x_{i-1}+x_{i+1}-2x_{i}\right)\:\forall x_{i}\in int\left(\Delta\right),\\
\left|g_{0}-u'\left(x_{0}\right)\right| & =a\left(x_{1}-x_{0}\right),\,\mathrm{and}\\
\left|g_{n}-u'\left(x_{n}\right)\right| & =a\left(x_{n-1}-x_{n}\right).
\end{alignat*}
\end{corollary}
\begin{remark}
For $i \in B_n$ let $h_i = x_{i+1} -x_i$. 
Then we have 
\[ \left|g_{i}-u'\left(x_{i}\right)\right|  =a\left(x_{i-1}+x_{i+1}-2x_{i}\right) = 
a (h_{i+1}-h_i). \]
We still get superapproximation of the gradient recovery when $h_{i+1}-h_i = O(h^2)$ 
when $i \in B_n$. 
\end{remark}

\subsection{Application to cubic functions}
\begin{corollary}
Let $u\in P_{3}(\mathbb{R})$ with $u(x)=ax^{3}+bx^{2}+cx+d$. Then,
\[
\left|g_{i}-u'\left(\tilde{x}_{i}\right)\right|=\frac{a}{4}\left(x_{i-1}-x_{i+1}\right)^{2}
\]
for all $i\in B_{n}$, and 
\begin{alignat*}{1}
\left|g_{0}-u'\left(\tilde{x}_{0}\right)\right| & =\frac{a}{4}\left(x_{0}-x_{1}\right)^{2}\\
\\
\left|g_{n}-u'\left(\tilde{x}_{n}\right)\right| & =\frac{a}{4}\left(x_{n-1}-x_{n}\right)^{2}
\end{alignat*}
where $\tilde{x}_{i}$ is defined as in Theorem \ref{th1}.
Similarly, for all $ i \in B_n$ we have 
\[ \left|g_{i}-u'(x_{i})\right|=a\left(x_{i-1}^{2}+x_{i-1}x_{i+1}+x_{i+1}^{2}-3x_{i}^{2}\right)+b\left(x_{i-1}+x_{i+1}-2x_{i}\right),
\] and 
\begin{align*} \left|g_{0}-u'(x_{0})\right|&=&a\left(x_{1}^{2}+x_{0}x_{1}-2x_{0}^{2}\right)+b\left(x_{1}-x_{0}\right), \\
\left|g_{n}-u'(x_{n})\right|&=&a\left(x_{n-1}^{2}+x_{n-1}x_{n}-2x_{n}^{2}\right)+b\left(x_{n-1}-x_{n}\right).\end{align*}
\end{corollary}
\Proof
The proof of this theorem is similar to Theorem \ref{th1}. 
\EndProof

\section{Conclusion}
We have presented an analysis of approximation property of the reconstructed gradient using 
an oblique projection. The reconstruction of the gradient is numerically efficient due to the use 
of a biorthogonal system. It is useful to investigate the extension  to higher order finite elements.

\bibliographystyle{plain}
\bibliography{total.bib}
\end{document}